\documentclass[11pt,a4paper]{article}
\usepackage[latin1]{inputenc}
\usepackage{latexsym}
\usepackage{amsmath,amssymb,amsthm}
\usepackage{rotating}
\usepackage{fancyhdr}
\usepackage{lscape}
\author{Barbara McGillivray\\University of Firenze, Italy}

\title{A probabilistic algorithm for the secant defect of Grassmann
varieties}
\date{}
\allowdisplaybreaks \textwidth = 13.3cm 
\newtheorem{teo}{Theorem}[section]
\newtheorem{cor}[teo]{Corollary}

\newtheorem{lemma}[teo]{Lemma}

\theoremstyle{definition}
\newtheorem{definiz}[teo]{Definition}

\begin{document}
\maketitle \pagenumbering{arabic} \baselineskip = 20pt

\begin{small}
\section{Abstract}

In this paper we study the higher secant varieties of Grassmann
varieties in relation to Waring's problem for alternating tensors
and to Alexander-Hirschowitz theorem. We show how to identify
defective higher secant varieties of Grassmannians using a
probabilistic method involving Terracini's Lemma, and we describe an
algorithm which can compute, by numerical methods, $dim(G(k,n)^{s})$
for $n\leq 14$. Our main result is that, except for Grassmannians of
lines, if $n\leq 14$ and $k\leq\frac{n-1}{2}$ (if $n=14$ we have
studied the case $k\leq 5$) there are only the four known defective
cases: $G(2,6)^{3}$, $G(3,7)^{3}$, $G(3,7)^{4}$ and $G(2,8)^{4}$.
\end{small}

\section{Introduction}
\emph{Waring's problem} for alternating tensors can be expressed in
the following form (see \cite{C})\\
\begin{quote}
Given a vector space $V$ of dimension $n+1$ and an alternating
tensor $\omega\in \bigwedge^{k+1}V$, what is the least integer $s$
such that $\omega$ can be written as the sum of $s$ decomposable
tensors of the form $v_{1}\wedge\ldots\wedge v_{k+1}$?
\end{quote}
This problem is still open and in this paper we will give some
evidence for what we expect the correct answer to be.\\
In order to formulate our result, we will consider a vector space
$V$ of dimension $n+1$ defined over a field $\mathbb{K}$ of
characteristic zero, and the Grassmann variety
$G(\mathbb{K}^{k+1},V)=G(k,n)$, which parametrises the decomposable
tensors in the projective space of $\bigwedge^{k+1}V$. As will be
explained in the next section, the problem translates into finding
the dimension of the $s$-secant variety $G(k,n)^{s}$ (see definition
\ref{definiz:secant}). The expected dimension of $G(k,n)^{s}$ is
$min\{\binom{n+1}{k+1}-1,s(n-k)(k+1)+s-1\}$, otherwise $G(k,n)^{s}$
is called \textbf{defective} (see definition
\ref{definiz:defective}).\\
It is well known that the Grassmannians of lines
$G(\mathbb{P}^{1},\mathbb{P}^{n})^{s}$ are defective until they fill
the ambient space and a list of four defective $G(k,n)^{s}$ is given
in \cite{CGG}. We would like to know if there exist other defective
varieties which are still unknown.\\
Computing the dimension of $G(k,n)^{s}$ is quite difficult, even
with the aid of a symbolic computation package; indeed just after
the defective examples of \cite{CGG}, the computer's memory reaches
its limit with the usual elimination technique using Gr\"{o}bner basis.\\
The main idea behind this paper is that one can compute
$dim(G(k,n)^{s})$ by means of a probabilistic method, which consists
in studying the span of the tangent spaces at $s$ chosen random
points. The dimension of this span can be computed by numerical
methods as the rank of a large matrix, and when this dimension
coincides with that expected, we can be sure that $G(k,n)^{s}$ is
not defective, indeed with another choice of points the dimension
cannot be larger because of inequality (\ref{eq:dif}).\\ This
technique allows us to take the computations further and our main
result is the following.
\begin{teo}\label{teo:tuttedifettive}
If $n\leq 14$ and $k\leq\frac{n-1}{2}$ \footnote{If we choose a
basis for $V$ there is a natural 1-1 correspondence between the
associated bases of $\bigwedge^{k+1}V$ and $\bigwedge^{n-k}V$, and
between the varieties $G(k,n)$ and $G(n-k-1,n)$. Thus we will only
consider the variety $G(k,n)$ where $k\leq\frac{n-1}{2}$.}(if $n=14$
we consider $k\leq 5$), $G(k,n)^{s}$ is defective only for
\begin{itemize}
\item $k=1$
\item $(k,n,s)=(2,6,3)$, $\delta=1$
\item $(k,n,s)=(3,7,3)$, $\delta=1$
\item $(k,n,s)=(3,7,4)$, $\delta=4$
\item $(k,n,s)=(2,8,4)$, $\delta=2$.
\end{itemize}
\end{teo}
This theorem is equivalent to the following answer to Waring's
problem:
\begin{cor}Given a finite dimension vector space $V$ and an
alternating tensor $\omega\in \bigwedge^{k+1}V$, if $n\leq 14$ and
$k\leq\frac{n-1}{2}$ (if $n=14$ we consider $k\leq 5$), then
$\omega$ can be written as the sum of $s$ decomposable tensors of
the form $v_{1}\wedge\ldots\wedge v_{k+1}$, where
$s=\lceil\frac{1}{(k+1)(n-k)+1}\binom{n+1}{k+1}\rceil$, except for
$(k,n,s)=(2,6,3)\,,(3,7,3)\,,(3,7,4)\,,(2,8,4)$ and $k=1$.
\end{cor}

\section{Waring's problem and some notations}

\emph{Waring's polynomial problem} has attracted considerable
attention from geometers and algebraists throughout its long and
absorbing history since it was first put forward in 1770. This
problem is connected with crucial issues in both
re-presentation theory and coding theory.\\
It poses the following question: if $f$ is a homogeneous polynomial
of degree $k$ in $n$ variables, what is the least integer $s$ such
that $f$ can be written as the sum of $k$th-powers
of $s$ linear forms?\\
This formulation of Waring's problem was solved in 1995 by J.
Alexander and A. Hirschowitz, who produced a formula for finding the
integer $s$; nevertheless this formula has four well known
exceptions.
\begin{teo}\label{teo:waring}\cite{AH}
Let $char(\mathbb{K})=0$. A homogeneous polynomial
$f\in\mathbb{K}[X_{0},\ldots,X_{n}]$ of degree $k$ can be
represented as the sum of $s$ powers of linear forms
\begin{displaymath}
f=L_{1}^{k}+\ldots+L_{s}^{k},
\end{displaymath}
where $s=\lceil\frac{1}{n+1}\binom{k+n}{n}\rceil$, except in the
cases where $(k,n,s)=(4,2,5)\,,\,(4,3,9)\,,\\(4,4,14)$\,,\,
$(3,4,7)$ and $k=2$.
\end{teo}
This challenging result can also be expressed in geometrical
terms, as we will explain.\\
First let us recall some fundamental definitions.\\
\newline
Let $X\subseteq\mathbb{P}^{N}$ be an $n$-dimensional irreducible
projective variety,
\begin{definiz}\label{definiz:secant}
The \textbf{$s$-secant variety} is the closure of the union of all
linear spaces spanned by $s$ points of $X$, which is expressed as
follows:\\
\begin{displaymath}
X^{s}=\overline{\bigcup_{x_{1},\ldots,x_{s}\in{X}}\,<x_{1},\ldots,x_{s}\,>}.
\end{displaymath}
\end{definiz}
The choice of s points in X gives rise to $sn$ free parameters. In
addition, $s$ points span a space of projective dimension $s-1$ and
$X^{s}$ will be embedded in $\mathbb{P}^{N}$. Consequently, we
should expect the dimension of $X^{s}$ to be given by
$min\{N,sn+s-1\}$: this is called the \textbf{expected
dimension }for secant varieties.\\
\newline
The following estimate on the dimension of $X^{s}$ is valid in
general:
\begin{equation}\label{eq:dif}
dim(X^{s})\le\min\{N,sn+s-1\}.
\end{equation}
From our viewpoint, the cases where the strict inequality applies
are the most interesting.
\begin{definiz}\label{definiz:defective}
The secant variety $X^{s}$ is called \textbf{defective} if
\begin{displaymath}
dim(X^{s})<min\{N,sn+s-1\}
\end{displaymath}
and the quantity $\delta=min\{N,sn+s-1\}-dim(X^{s})$ is its
\textbf{defectiveness}.
\end{definiz}
One well-known result that is useful in finding the dimensions of
the multisecant varieties is \emph{Terracini's lemma}.
\begin{lemma}[\cite{Z}]
Let $x_{1},\ldots,x_{s}\in X$ be generic points; let us refer to the
projectivised tangent spaces to $X$ at these points as
$T_{x_{1}}X,\ldots,$ $T_{x_{s}}X$, then
\begin{displaymath}
dim(X^{s})=dim(<T_{x_{1}}X,\ldots,T_{x_{s}}X>).
\end{displaymath}
\end{lemma}
Let us now take a homogeneous polynomial $f(X_{0},\ldots,X_{n})$ of
degree $k$.\\ Asking whether $f$ can be written as the sum of powers
of degree $k$ of $s$ linear forms $L_{1},\ldots,L_{s}$ is the same
as asking whether $f$ belongs to the $s$-secant variety of the $k$th
Veronese embedding of $\mathbb{P}^{n}$, which we call $V_{k,n+1}$.
It is therefore important to know the dimension of $V_{k,n+1}^{s}$,
and consequently the cases where $V_{k,n+1}$ is defective for
$s$-secant varieties.\\
The result obtained by Alexander and Hirschowitz is extremely useful
in our case and translates geometrically as follows:
\begin{teo}(\cite{IK})
The Veronese variety $V_{k,n+1}$ is defective for $s$-secant
varie-ties only in the following cases:
\begin{displaymath}
(k,n,s)=(4,2,5)\,,\,(4,3,9)\,,\,(4,4,14)\,,\,(3,4,7)
\end{displaymath}
and $k=2$.
\end{teo}
We have therefore obtained a full classification of defective
Veronese varieties.\\
\newline
In this paper we will analyse the problem of defectiveness with
respect to another important family of classical varieties, the
Grassmannians,\footnote{For the problem of defectiveness of Segre
Varieties and its connection with the rank of tensors, see
\cite{S}.}
which are related to exterior algebras.\\
If $k\leq \dim(V)$ is a positive integer, we define the
\textbf{Grassmannian} $G(k,V)$ to be the variety of projective
subspaces of $\mathbb{P}(V)$ of dimension $k$. When
$V=\mathbb{K}^{n+1}$, $G(k,V)$ will be denoted by
$G(k,n)$.\\
Since $dim(G(k,n))=(k+1)(n-k)$, the expected dimension for the
secant varieties $G(k,n)^{s}$ is:
\begin{displaymath}
min\Big\{\binom{n+1}{k+1}-1,s(n-k)(k+1)+s-1\Big\}
\end{displaymath}
We also have the two following important theoretical results to draw
on.
\begin{teo}\label{teo:k=1}
$G(1,n)^{s}$ is defective for $s<\lfloor{\frac{n}{2}}\rfloor$.
\end{teo}
\begin{teo} (\cite{CGG}) \label{teo:casi}
Let $k\geq 2$. If $s(k+1)\leq n+1$ , then $G(k,n)^{s}$ has the
expected dimension.
\end{teo}

\section{A probabilistic algorithm and proof of Theorem\ref{teo:tuttedifettive}}

To tackle our problem, we initially used the \emph{Macaulay 2}
computation system (see \cite{MAC}), which was designed to study
problems of algebraic geometry and commutative algebra, to write an
algorithm generating parametric equations for the Grassmannians we
are studying (see \cite{B}).\\
To calculate the dimensions of the multisecant varieties, we
favoured a proba-bilistic approach involving \emph{Terracini's
lemma}. We took $s$ random points in $G(k,n)$ and studied their
tangent spaces and the space spanned by these tangent spaces. If we
found the expected dimension, the result was clearly correct, but if
this revealed defectiveness, more checks needed to be performed.\\
Using this approach we constructed an algorithm that turned our
problem into the calculation of the rank of fairly large matrices
with constant coefficients; to study $dim(G(k,n)^{s})$ we needed to
know the rank of a matrix of order $s(1+(k+1)(n-k))\times N$. This
algorithm enabled us to compute the dimension of $G(k,n)^{s}$ when
$n\leq 11$, $k\leq 4$ and $s\leq 4$, at which stage the computer's
memory was used up. It was therefore clear that symbolic computation
was not the best tool for this type of task.\\
To further proceed with our study, we decided to employ the
\emph{Matlab} software system, which is a computation system
designed for dealing with numerical computations involving very large matrices.\\
The new algorithm obtained confirms the validity of the
probabilistic approach. It is based on theoretical observation that
the tangent spaces can be computed without having to define
equations for the Grassmannian.\\
Let us take a point $P=v_{0}\wedge\ldots\wedge v_{k}\in G=G(k,n)$;
the Plucker coordinates of $P$ are all the $(k+1)\times(k+1)$ minors
of the matrix $A$ of order $(k+1)\times(n+1)$, which has the vectors
$v_{0},\ldots,v_{k}\in V$ for rows. It is easy to check, by the
Leibniz rule, that the following is true:
\begin{lemma}
$T_{P}(G)$ is the projective space associated with
\begin{displaymath}
V\wedge v_{1}\wedge\ldots\wedge v_{k}\,+\,v_{0}\wedge V\wedge
v_{2}\wedge\ldots\wedge v_{k}\,+\,\ldots\,+\,v_{0}\wedge\ldots\wedge
v_{k-1}\wedge V\, =\,T_{0}\,+\,\ldots\,+\,T_{k}.
\end{displaymath}
\end{lemma}
If $A_{i,j}$ stands for the matrix obtained from $A$ by replacing
the $i$th row by the $j$th row of the identity matrix $I$ of order
$(n+1)\times(n+1)$, then every $T_{i}$ is parametrised by a matrix
$M_{i}$ of order $(n+1)\times N$ whose $j$th row $m_{j}$ contains
the minors of maximum order of the matrix $A_{i,j}$.\\
\begin{small}
\begin{equation}\label{eqn:aij}
\begin{array}{ccc}A_{i,1}\,=\,\left ( \begin{array}{cccc}
&v_{0}&\\&\vdots&&\\&v_{i-1}&\\1&0&\ldots&0\\
&v_{i+1}&\\&\vdots&&\\&v_{k}&
\end{array}\right),\,
\ldots\,,\,& \left ( \begin{array}{ccccc}
&v_{0}&\\&\vdots&&\\&v_{i-1}&\\0&\ldots&0&1\\
&v_{i+1}&\\&\vdots&&\\&v_{k}&
\end{array}\right)=A_{i,n+1}
\end{array}
\end{equation}
\begin{displaymath}
M_{i}\,=\,\left ( \begin{array}{c} m_{1}\\
\vdots\\ m_{n+1}
\end{array}\right)
\end{displaymath}
\end{small}
Our algorithm is described below.

\begin{itemize}
\item Input: positive integers $n_{m}$ and $n_{M}$.
\item Repeat on parameters $n_{m}\leq n\leq n_{M}$, $1\leq k\leq
\lfloor\frac{n-1}{2}\rfloor$ and $2\leq s\leq
S=\lceil\frac{1}{(k+1)(n-k)+1}\binom{n+1}{k+1}\rceil$) to study
$G^{s}$.
\item Define the matrix $TA$ that contains the actual dimensions
and the defectiveness of $G^{s}$.
\item Define the function $ed(k,n,s)$ that
calculates the expected dimension of $G^{s}$ and the matrix $E$ of
the expected dimensions.
\item Choose $s$ random points $P_{1},\ldots,P_{s}$ in $G$.
\begin{itemize}
\item Take a matrix $B$ of order $s(k+1)\times(n+1)$ with random rational coefficients
in the interval $[-L,L]$.
\item Extract $s$ submatrices $A$ of order $(k+1)\times(n+1)$ from $B$.
\end{itemize}
\item Repeat for $1\leq h\leq s$ and study
$T_{P_{h}}G=T_{1}+\ldots+T_{k+1}$.
\begin{itemize}
\item Repeat for $1\leq i\leq k+1$.
\item For every $1\leq j\leq n+1$ calculate the minors $(k+1)\times(k+1)$
of $A_{i,j}$, computed from $A$ as in (\ref{eqn:aij}), and call the
row of these minors $m_{j}$.
\item
Construct the matrix $M_{i}$ with rows $m_{j}$.
\end{itemize}
\item Parametrise $T_{P_{1}}G+\ldots+T_{P_{s}}G$.
\begin{itemize}
\item Concatenate $M_{1},\ldots,M_{k+1}$ vertically to
obtain matrix $M$.
\end{itemize}
\item Determine the value of the projective dimension of $G^{s}$.
\begin{itemize}
\item Calculate the rank of $M$, then subtract 1.
\item Define row $dim$ of $TA$ of actual dimensions
and row $dif=E-dim$ of defectiveness.
\end{itemize}
\item Output: matrix $TA$.
\end{itemize}
This is the text of the algorithm. \linespread{1.17}
\begin{footnotesize}
\begin{verbatim}

L=100
nm=3
nM=14
f=inline('floor((n-1)/2)')
F=[]
for u=nm:nM
F=[F 2*f(u)]
end
N=inline('factorial(n+1)/(factorial(k+1)*factorial(n-k))-1')
S=inline('ceil((N+1)/((n-k)*(k+1)+1))')
Smax=S(N(f(nM),nM),f(nM),nM)
TA=[]
  for n=nm:nM
  T=zeros(2*f(n),Smax+1)
    for k=1:f(n)
    E=[]
    T(2*k-1,1:2)=[k n]
    I=eye(n+1)
    v=nchoosek(1:n+1,k+1)
    l=size(v,1)
    dim=[]
      for s=2:S(N(k,n),k,n)
      M=[]
      ed=inline('min(N,s*(n-k)*(k+1)+s-1)')
      E=[E ed(N(k,n),k,n,s)]
      B=rand(k+1,(n+1)*s)
      B=(B-0.5)*2*L
        for h=1:s
        A=B(:,(h-1)*(n+1)+1:h*(n+1))
          for i=1:k+1
          r=A(i,:)
            for j=1:n+1
            A(i,:)=I(j,:)
            m=[]
              for w=1:l
              D(w)=det(A(:,v(w,:)))
              m=[m D(w)]
              end
            M=[M;m]
            end
          A(i,:)=r
          end
        end
      dim=[dim rank(M)-1]
      dif=E-dim
      T(2*k-1,s+1)=dim(1,s-1)
      T(2*k,s+1)=dif(1,s-1)
      end
    end
    TA=[TA;T]
  end
TA

\end{verbatim}
\end{footnotesize}

If $s=S=\lceil\frac{1}{(k+1)(n-k)+1}\binom{n+1}{k+1}\rceil$ and
$G(k,n)^{s}$ is not defective, then the variety fills the ambient
space $\mathbb{P}^{N}=\mathbb{P}^{\binom{n+1}{k+1}-1}$; if it is
defective, then we find $dif>0$ and it can happen that
$dim(G(k,n)^{S})<N$, so that we will have to calculate $dim(G(k,n)^{s})$ for $s>S$.\\
Using this algorithm, at the stage $(n,k)=(6,14)$ the computer's
memory was used up. Our results are summarised in the following
tables.

\oddsidemargin = 1.48cm
\begin{center}
\begin{footnotesize}
\begin{tabular}{||c c |c c||c|c|c|c|c|c|c|c|c|c|c|c|c|c||}
\hline
$N$&$S$&$k$&$n$&$G^{2}$&$G^{3}$&$G^{4}$&$G^{5}$&$G^{6}$&$G^{7}$&$G^{8}$&$G^{9}$&$G^{10}$&
$G^{11}$&$G^{12}$&$G^{13}$&$G^{14}$&$G^{15}$
\\ \hline5&2&1&3&5&&&&&&&&&&&&&
\\ \hline9&2&1&4&9 &     &&&&&&&&&&&&\\
\hline14  &2&1&5&13*&14   &&&&&&&&&&&&\\
\hline19&2&2&5&19&&&&&&&&&&&&&\\
\hline20&2&1&6&17*&20&&&&&&&&&&&&\\
\hline34&3&2&6&25&33*&34&&&&&&&&&&&\\
\hline27&3&1&7&21*&26*&27&&&&&&&&&&&\\
\hline55&4&2&7&31&47&55&&&&&&&&&&&\\
\hline69&5&3&7&33&49*&63*&69&&&&&&&&&&\\
\hline35&3&1&8&25*&32*&35&&&&&&&&&&&\\
\hline83&5&2&8&37&56&73*&83&&&&&&&&&&\\
\hline125&6&3&8&41&62&83&104&125&&&&&&&&&\\
\hline44&3&1&9&29*&38*&43*&44&&&&&&&&&&\\
\hline119&6&2&9&43&65&87&109&119&&&&&&&&&\\
\hline209&9&3&9&49&74&99&124&149&174&199&209&&&&&&\\
\hline251&10&4&9&51&77&103&129&155&181&207&233&251&&&&&\\
\hline54&3&1&10&33*&44*&51*&54&&&&&&&&&&\\
\hline164&7&2&10&49&74&99&124&149&164 &  &&&&&&&\\
\hline329&12&3&10&57&86&115&144&173&202&231&260&289&318&329&&&\\
\hline461&15&4&10&61&92&123&154&185&216&247&278&309&340&371&402&433&461\\
\hline65&4&1&11&37*&50*&59*&64*&65&&&&&&&&&\\
\hline219&8&2&11&55&83&111&139&167&195&219&&&&&&&\\
\hline494&15&3&11&65&98&131&164&197&230&263&296&329&362&395&428&461&494\\
\hline791&22&4&11&71&107&143&179&215&251&287&323&359&395&431&467&503&539\\
\hline923&25&5&11&73&110&147&184&221&258&295&332&369&406&443&480&517&554\\
\hline77 &6 &1&12&41*&56* &67* &74* &77 &   &   &   &   &   &   &   &   &\\
\hline285&10&2&12&61&92 &123&154&185&216&247&278&285&   &   &   &   &\\
\hline714&20&3&12&73&110&147&184&221&258&295&332&369&406&443&480&517&554\\
\hline1286&32&4&12&81&122&163&204&245&286&327&368&409&450&491&532&573&614\\
\hline1715&40&5&12&85&128&171&214&257&300&343&386&429&472&515&558&601&644\\
\hline90  &4 &1&13&45*&62* &75* &84* &89* &90 &   &   &   &   &   &   &   &\\
\hline363&11&2&13&67&101&135&169&203&237&271&305&339&363&&&&\\
\hline1000&25&3&13&81&122&163&204&245&286&327&368&409&450&491&532&573&614\\
\hline2001&44&4&13&91&137&183&229&275&321&367&413&459&505&551&597&643&689\\
\hline3002&62&5&13&97&146&195&244&293&342&391&440&489&538&587&636&685&734\\
\hline3431&69&6&13&99&149&199&249&299&349&399&449&499&549&599&649&699&749\\
\hline104&4&1&14&49*&68*&83*&94*&101*&104&&&&&&&&\\
\hline454&13&2&14&73&110&147&184&221&258&295&332&369&406&443&454&&\\
\hline1364&31&3&14&89 &134&179&224&269&314&359&404&449&494&539&584&629&674\\
\hline3002&59&4&14&101&152&203&254&305&356&407&458&509&560&611&662&713&764\\
\hline5004&91&5&14&109&164&219&274&329&384&439&494&549&604&659&714&769&824\\
\hline
\end{tabular}
\end{footnotesize}
\end{center}

\newpage
\begin{center}
\begin{footnotesize}
\begin{tabular}{||c c |c c||c|c|c|c|c|c|c|c|c|c|c|c|c||}
\hline
$N$&$S$&$k$&$n$&$G^{16}$&$G^{17}$&$G^{18}$&$G^{19}$&$G^{20}$&$G^{21}$
&$G^{22}$&$G^{23}$&$G^{24}$&$G^{25}$ &$G^{26}$&$G^{27}$&$G^{28}$\\
\hline791&22&4&11&575&611&647&683&719&755&791&&&&&&\\
\hline923&25&5&11&591&628&665&702&739&776&813&850&887&923&&&\\
\hline714&20&3&12&591&628&665&702&714&&&&&&&&\\
\hline1286&32&4&12&655&696&737&778&819&860&901&942&983&1024&1065&1106&1147\\
\hline1715&40&5&12&687&730&773&816&859&902&945&988&1031&1074&1117&1160&1203\\
\hline1000&25&3&13&655&696&737&778&819&860&901 &942 &983 &1000&&\\
\hline2001&44&4&13&735&781&827&873&919&965&1011&1057&1103&1149&1195&1241&1287\\
\hline3002&62&5&13&783&832&881&930&979&1028&1077&1126&1175&1224&1273&1322&1371\\
\hline3431&69&6&13&799&849&899&949&999&1049&1099&1149&1199&1249&1299&1349&1399\\
\hline1364&31&3&14&719&764&809&854&899&944&989&1034&1079&1124&1169&1214&1259\\
\hline3002&59&4&14&815&866&917&968&1019&1070&1121&1172&1223&1274&1325&1376&1427\\
\hline5004&91&5&14&879&934&989&1044&1099&1154&1209&1264&1319&1374&1429&1484&1539\\
\hline
\end{tabular}

\begin{tabular}{||c c |c c||c|c|c|c|c|c|c|c|c|c|c|c||}
\hline
$N$&$S$&$k$&$n$&$G^{29}$&$G^{30}$&$G^{31}$&$G^{32}$&$G^{33}$&$G^{34}$
&$G^{35}$&$G^{36}$&$G^{37}$&$G^{38}$ &$G^{39}$&$G^{40}$\\
\hline1286&32&4&12&1188&1229&1270&1286&    &    &    &    &    &    &    &\\
\hline1715&40&5&12&1246&1289&1332&1375&1418&1461&1504&1547&1590&1633&1676&1715\\
\hline2001&44&4&13&1333&1379&1425&1471&1517&1563&1609&1655&1701&1747&1793&1839\\
\hline3002&62&5&13&1420&1469&1518&1567&1616&1665&1714&1763&1812&1861&1910&1959\\
\hline3431&69&6&13&1449&1499&1549&1599&1649&1699&1749&1799&1849&1899&1949&1999\\
\hline1364&31&3&14&1304&1349&1364&&&&&&&&&\\
\hline3002&59&4&14&1478&1529&1580&1631&1682&1733&1784&1835&1886&1937&1988&2039\\
\hline5004&91&5&14&1594&1649&1704&1759&1814&1869&1924&1979&2034&2089&2144&2199\\
\hline
\end{tabular}

\begin{tabular}{||c c |c c||c|c|c|c|c|c|c|c|c|c|c|c||}
\hline
$N$&$S$&$k$&$n$&$G^{41}$&$G^{42}$&$G^{43}$&$G^{44}$&$G^{45}$&$G^{46}$
&$G^{47}$&$G^{48}$&$G^{49}$&$G^{50}$ &$G^{51}$&$G^{52}$\\
\hline2001&44&4&13&1885&1931&1977&2001&&&&&&&&\\
\hline3002&62&5&13&2008&2057&2106&2155&2204&2253&2302&2351&2400&2449&2498&2547\\
\hline3431&69&6&13&2049&2099&2149&2199&2249&2299&2349&2399&2449&2499&2549&2599\\
\hline3002&59&4&14&2090&2141&2192&2243&2294&2345&2396&2447&2498&2549&2600&2651\\
\hline5004&91&5&14&2254&2309&2364&2419&2474&2529&2584&2639&2694&2749&2804&2859\\
\hline
\end{tabular}

\begin{tabular}{||c c |c c||c|c|c|c|c|c|c|c|c|c|c|c||}
\hline
$N$&$S$&$k$&$n$&$G^{53}$&$G^{54}$&$G^{55}$&$G^{56}$&$G^{57}$&$G^{58}$
&$G^{59}$&$G^{60}$&$G^{61}$&$G^{62}$ &$G^{63}$&$G^{64}$\\
\hline3002&62&5&13&2596&2645&2694&2743&2792&2841&2890&2939&2988&3002&&\\
\hline3431&69&6&13&2649&2699&2749&2799&2849&2899&2949&2999&3049&3099&3149&3199\\
\hline3002&59&4&14&2702&2753&2804&2855&2906&2957&3002&    &    &    &    &\\
\hline5004&91&5&14&2914&2969&3024&3079&3134&3189&3244&3299&3354&3409&3464&3519\\
\hline
\end{tabular}

\begin{tabular}{||c c |c c||c|c|c|c|c|c|c|c|c|c|c|c||}
\hline
$N$&$S$&$k$&$n$&$G^{65}$&$G^{66}$&$G^{67}$&$G^{68}$&$G^{69}$&$G^{70}$
&$G^{71}$&$G^{72}$&$G^{73}$&$G^{74}$ &$G^{75}$&$G^{76}$\\
\hline3431&69&6&13&3249&3299&3349&3399&3431&    &    &    &    &    &    &\\
\hline5004&91&5&14&3574&3629&3684&3739&3794&3849&3904&3959&4014&4069&4124&4179\\
\hline
\end{tabular}

\begin{tabular}{||c c |c c||c|c|c|c|c|c|c|c|c|c|c|c||}
\hline
$N$&$S$&$k$&$n$&$G^{77}$&$G^{78}$&$G^{79}$&$G^{80}$&$G^{81}$&$G^{82}$
&$G^{83}$&$G^{84}$&$G^{85}$&$G^{86}$ &$G^{87}$&$G^{88}$\\
\hline5004&91&5&14&4234&4289&4344&4399&4454&4509&4564&4619&4674&4729&4784&4839\\
\hline
\end{tabular}

\begin{tabular}{||c c |c c||c|c|c||}
\hline $N$&$S$&$k$&$n$&$G^{89}$&$G^{90}$&$G^{91}$\\
\hline5004&91&5&14&4894&4949&5004\\
\hline
\end{tabular}
\end{footnotesize}
\end{center}

\newpage
\begin{normalsize}
These results confirm what is known about the defectiveness of
Grassmannians. Except for Grassmannians of lines, we have identified
four defective varieties. We have therefore proved theorem
\ref{teo:tuttedifettive}.
\end{normalsize}


\addcontentsline{toc}{section}{References}
\begin{thebibliography}{99}
\bibitem{C}C. Ciliberto, \emph{Geometric aspects of polynomial interpolation in more variables and
of Waring's problem}, Proceedings of the European Congress in
Mathematics, Vol. I (Barcelona, 2000), 289-316, Birkh\"{a}user 2001
\bibitem{CGG}M. V. Catalisano, A. V. Geramita, A. Gimigliano,
\emph{Secant varieties of Grassmann Varieties}, Proceedings of the
American Mathematical Society (math.AG/0208166).
\bibitem{AH}J. Alexander, A. Hirschowitz, \emph{Polynomial
interpolation in several variables}, Journal of Alg. Geom.
\textbf{4} (1995), 201-222.
\bibitem{Z}F. L. Zak, \emph{Tangents and secants of algebraic varieties},
Translations of \\Mathematical Monographs, Vol. 127, American
Mathematical Society, 1993.
\bibitem{IK}A. Iarrobino, V. Kanev, Power Sums, \emph{Gorenstein Algebras, and
Determinantal Loci}, Springer-Verlag, Berlin-Heidelberg, 1999.
\bibitem{S}V. Strassen, \emph{Rank and optimal computation of generic tensors},
Linear\\
Algebra and its Applications \textbf{52/53} (1983), 645-685.
\bibitem{MAC}D. R. Grayson, M. E. Stillman, \emph{Macaulay 2}, Software system available on website www.math.uiuc.edu.
\bibitem{B}B. McGillivray, \emph{Metodi computazionali per lo studio delle variet\`{a}
di Grassmann}, Laurea Thesis, Florence, 2004.


\end{thebibliography}
\end{document}